\documentclass[12pt,reqno]{amsart}
\usepackage{amssymb}
\usepackage{amsmath}
\usepackage{amsthm}
\usepackage{eucal}
\usepackage{color}
\usepackage{amsfonts}
\usepackage[pdftex]{graphicx}
\usepackage[T1]{fontenc}

\newcommand{\car}{\chi_{_{\epsilon,b}}}
\newcommand{\R}{\mathbb R}

\newcommand{\Z}{\mathbb Z}

\newtheorem{theorem}{Theorem}[section]
\newtheorem{proposition}[theorem]{Proposition}
\newtheorem{remark}[theorem]{Remark}

\newtheorem{lemma}[theorem]{Lemma}
\newtheorem{corollary}[theorem]{Corollary}

\newtheorem*{TA1}{Theorem A1}
\newtheorem*{TA2}{Theorem A2}
\newtheorem*{TA3}{Theorem A3}
\newtheorem*{TA4}{Theorem A4}
\newtheorem*{L2}{Lemma 2}

\numberwithin{equation}{section}
\begin{document}
\title[Regularity of Solutions]{On the regularity of  solutions to the $k$-generalized Korteweg-de Vries equation}
\author{C. E. Kenig}
\address[C. E. Kenig]{Department of Mathematics\\University of Chicago\\Chicago, Il. 60637 \\USA.}
\email{cek@math.uchicago.edu}
\author{F. Linares}
\address[F. Linares]{IMPA\\
Instituto Matem\'atica Pura e Aplicada\\
Estrada Dona Castorina 110\\
22460-320, Rio de Janeiro, RJ\\Brazil}
\email{linares@impa.br}

\author{G. Ponce}
\address[G. Ponce]{Department  of Mathematics\\
University of California\\
Santa Barbara, CA 93106\\
USA.}
\email{ponce@math.ucsb.edu}

\author{L. Vega}
\address[L. Vega]{UPV/EHU\\Dpto. de Matem\'aticas\\Apto. 644, 48080 Bilbao, Spain, and Basque Center for Applied Mathematics,
E-48009 Bilbao, Spain.}
\email{luis.vega@ehu.es}


\begin{abstract} This work is concerned with special regularity properties of solutions to the $k$-generalized Korteweg-de Vries equation. In \cite{ILP-cpde} it was established that if the initial datun $u_0\in H^l((b,\infty))$ for some $l\in\Z^+$ and $b\in \R$, then  the corresponding solution 
$u(\cdot,t)$ belongs to $H^l((\beta,\infty))$ for any $\beta \in \R$ and any $t\in (0,T)$. Our
 goal here  is to extend this result to the case where $\,l>3/4$.
\end{abstract}

\maketitle

\section{Introduction}

In this note we study  the regularity of solutions to the initial value problem (IVP) associated to 
the $k$-generalized Korteweg-de Vries equation
\begin{equation}\label{AA}
\begin{cases}
\partial_tu+\partial_x^3u +u^k\partial_xu=0, \hskip5pt \;x, t\in\R, \; k\in\Z^{+},\\
u(x,0)=u_0(x).
\end{cases}
\end{equation}

The starting point is a property found by Isaza, Linares and Ponce \cite{ILP-cpde} concerning the  propagation of smoothness in solutions of the IVP \eqref{AA}.
  To state it we first recall the following well-posedness (WP) result for the IVP \eqref{AA}:
\begin{TA1}\label{A1}
If $u_0\in H^{{3/4}^{+}}(\R)$, then there exist $ T=T(\|u_0\|_{_{{\frac34}^{+}, 2}}; k)>0$ and a unique solution $u=u(x,t)$ of the IVP \eqref{AA}
such that
\begin{equation}\label{notes-2}
\begin{split}
{\rm (i)}\hskip25pt & u\in C([-T,T] : H^{{3/4}^{+}}(\R)),\\
{\rm(ii)}\hskip25pt & \partial_x u\in L^4([-T,T]: L^{\infty}(\R)), \hskip10pt {\text (Strichartz)},\\
{\rm(iii)}\hskip25pt &\underset{x}{\sup}\int_{-T}^{T} |J^r\partial_x u(x,t)|^2\,dt<\infty \text{\hskip10pt for \hskip10pt} r\in[0,{3/4}^{+}],\\
{\rm(iv)}\hskip25pt & \int_{-\infty}^{\infty}\; \sup_{-T\leq t\leq T}|u(x,t)|^2\,dx < \infty,
\end{split}
\end{equation}
with $\;J=(1-\partial_x^2)^{1/2}$. Moreover, the map data-solution, $\,u_0\to \,u(x,t)$ is locally continuous (smooth) 
from $\,H^{3/4+}(\R)$ into the class $\,X^{3/4+}_T\,$ defined in \eqref{notes-2}.

If $k\geq2$, then the result holds in $H^{{3/4}}(\R)$. If $k=1,2,3$, then $T$ can be taken arbitrarily large.
\end{TA1}

For the proof of Theorem A1 we refer to  \cite{KPV93}, \cite{CKSTT} and \cite{AG}. The proof of our main result Theorem \ref{A3} is based on an energy estimate argument for which the estimate (ii) in \eqref{notes-2} (i.e. the time integrability of $\|\partial_xu(\cdot,t)\|_{\infty}$)  is essential. 
However, we remark  that from the WP point of view is not optimal. For a detailed discussion on the WP of the IVP  \eqref{AA} we refer to \cite{LP}, Chapters 7-8.

Now we enunciate the result obtained in \cite{ILP-cpde} regarding propagation of regularities which motivates our study here:

\begin{TA2}[\cite{ILP-cpde}]\label{A2}
Let $u_0\in H^{{3/4}^{+}}(\R)$. If for some $\,l\in \Z^{+}$ and for some $x_0\in \R$
\begin{equation}\label{notes-3}
\|\,\partial_x^l u_0\|^2_{L^2((x_0,\infty))}=\int_{x_0}^{\infty}|\partial_x^l u_0(x)|^2dx<\infty,
\end{equation}
then the solution $u=u(x,t)$ of the IVP \eqref{AA} provided by Theorem A1 satisfies  that for any $v>0$ and $\epsilon>0$
\begin{equation}\label{notes-4}
\underset{0\le t\le T}{\sup}\;\int^{\infty}_{x_0+\epsilon -vt } (\partial_x^j u)^2(x,t)\,dx<c,
\end{equation}
for $j=0,1, \dots, l$ with $c = c(l; \|u_0\|_{{3/4}^{+},2};\|\,\partial_x^l u_0\|_{L^2((x_0,\infty))} ; v; \epsilon; T)$.

In particular, for all $t\in (0,T]$, the restriction of $u(\cdot, t)$ to any interval of the form $(a, \infty)$ belongs to $H^l((a,\infty))$.

Moreover, for any $v\geq 0$, $\epsilon>0$ and $R>0$ 
\begin{equation}\label{notes-5}
\int_0^T\int_{x_0+\epsilon -vt}^{x_0+R-vt}  (\partial_x^{l+1} u)^2(x,t)\,dx dt< c,
\end{equation}
with  $c = c(l; \|u_0\|_{_{{3/4}^{+},2}};\|\,\partial_x^l u_0\|_{L^2((x_0,\infty))} ; v; \epsilon; R; T)$.
\end{TA2}

Theorem A2  tells us that the $H^l$-regularity ($l\in\Z^+$) on the right hand side of the data travels forward in time with infinite speed. Notice that since the equation is reversible in time
a gain of regularity in $H^s(\R)$ cannot occur so at $t>0$, so  $u(\cdot, t)$ fails to be in $H^l(\R)$ due to its decay at $-\infty$. In
this regard, it was also shown  in  \cite{ILP-cpde} that for any  $\delta>0$ and $t\in (0,T)$ and $j=1,\dots,l$
\begin{equation*}
\int_{-\infty}^{\infty} \frac{1}{\langle x_{-}\rangle^{j+\delta}} (\partial_x^j u)^2(x,t)\,dx \le \frac{c}{t},
\end{equation*}
with $c= c(\|u_0\|_{{3/4}^{+},2}; \|\partial_x^j u_0\|_{L^2((x_0,\infty))}; \,x_0; \,\delta)$, $x_{-}=\max\{0;-x\}$ and $\langle x\rangle=(1+x^2)^{1/2}$.


 The aim of this note is to extend  Theorem A2 to the case where the local regularity of the datum $u_0$ in \eqref{notes-3}  is measure with a fractional exponent. Thus, our main result is :

\begin{theorem}\label{A3}
Let $u_0\in H^{{3/4}^{+}}(\R)$.  If for some $\,s\in\R,\,s>3/4,$ and for some $x_0\in \R$
\begin{equation}\label{anotes-3}
\|\,J^s u_0\|^2_{L^2((x_0,\infty))}=\int_{x_0}^{\infty}|J^s u_0(x)|^2dx<\infty,
\end{equation}
then the solution $\,u=u(x,t)\,$ of the IVP \eqref{AA} provided by Theorem A1 satisfies  that for any $v>0$ and $\epsilon>0$
\begin{equation}\label{anotes-4}
\underset{0\le t\le T}{\sup}\;\int^{\infty}_{x_0+\epsilon -vt } (J^ru)^2(x,t)\,dx<c,
\end{equation}
for $\,r\in(3/4,s]$ with $c = c(l; \|u_0\|_{{3/4}^{+},2};\|\,J^r u_0\|_{L^2((x_0,\infty))} ; v; \epsilon; T)$.

Moreover, for any $v\geq 0$, $\epsilon>0$ and $R>0$ 
\begin{equation}\label{anotes-5}
\int_0^T\int_{x_0+\epsilon -vt}^{x_0+R-vt}  (J^{s+1} u)^2(x,t)\,dx dt< c,
\end{equation}
with  $c = c(l; \|u_0\|_{_{{3/4}^{+},2}};\|\,J^s u_0\|_{L^2((x_0,\infty))} ; v; \epsilon; R; T)$.
\end{theorem}

From the results in section 2 it will be clear that we need only consider the case
$\,s\in (3/4,\infty)-\Z^{+}$.

The rest of this paper is organized as follows : section 2 contains some preliminary estimates required for Theorem \ref{A3}, whose proof will be given in section 3.

\section{Preliminary estimates}

Let $T_a$ be a pseudo-differential operator whose symbol 
\begin{equation}
\label{symbol}
\sigma(T_a)=a(x,\xi)\in S^r,\;\;r\in\R,
\end{equation}
 so that
\begin{equation}
\label{pseudo}
T_af(x)=\int_{\R^n}\,a(x,\xi)\widehat{f}(\xi)\,e^{2\pi i x\cdot \xi}d\xi.
\end{equation}
\vskip.1in

The following result is the singular integral realization of a pseudo-differential operator, whose proof can be found in \cite{St2} Chapter 4.

\begin{TA3}
\label{TA3}
Using the above notation \eqref{symbol}-\eqref{pseudo} one has that 
\begin{equation}
\label{pseudo-si}
T_af(x)=\int_{\R^n}\,k(x,x-y) f(y) dy,\;\;\;\;\;\text{if}\;\;\;\;x\notin supp (f)
\end{equation}
where $k\in C^{\infty}(\R^n\times R^n-\{0\})$ satisfies : $\forall \,\alpha,\beta \in (Z^+)^n\, \,\forall\,N\geq 0$ 
\begin{equation}
\label{estimate1}
|\partial_x^{\alpha}\partial_z^{\beta} k(x,z)|\leq A_{\alpha,\beta,N, \delta}\,|z|^{-(n+m+|\beta|+N)},\;\;\;\;\;\;\;\;\;\;|z|\geq \delta,
\end{equation}
\begin{equation*}
\;\text{if}\;\;\;\;\;\;n+m+|\beta|+N>0\;\;\;\;\;\;\;\;\;\;\;\text{uniformly in }\;\;x\in\R^n.
\end{equation*}
\end{TA3}
\vskip.1in
 To simplify the exposition from now on we restrict ourselves to the one-dimensional case $\,x\in\R$ where in the next section these results will be applied.
 
 As a direct consequence of Theorem A3 one has :
 \vskip.1in
\begin{corollary} 
\label{c1}
Let $m\in\Z^{+}$ and $l\in\R$. If $\,g\in L^2(\R)$ and $\,f\in L^p(\R)$, $p\in[2,\infty]$, with 
$$
distance( supp(f);\,supp(g))\geq \delta>0,
$$
then
\begin{equation}
\label{abc}
\| f\,\partial_x^m J^lg\|_2\leq c\|f\|_p\|g\|_2.
\end{equation}
\end{corollary}

\vskip.1in

Next, let $\,\theta_j\in C^{\infty}(\R)-\{0\}$ with $\,\theta_j'\in C_0^{\infty}(\R)$ for $\,j=1,2$ and 
\begin{equation}
\label{teta}
\,distance (supp(1-\theta_1);\,supp(\theta_2))\geq \delta>0.
\end{equation}
\vskip.1in

\begin{lemma}
\label {L1}
Let $\,f\in H^s(\R),\,s<0,$ and $\, T_a$ be a pseudo-differential operator of order zero ($a\in S^0$). If $\,\theta_1f\in L^2(\R)$, then 
\begin{equation}
\label{11}
\,\theta_2T_a f\in L^2(\R).
\end{equation}

\end{lemma}

\noindent{\bf Proof of Lemma \ref{L1}.} Since
\begin{equation*}
\theta_2\,T_af= \theta_2 T_a(\theta_1f) + \theta_2 T_a((1-\theta_1)f),
\end{equation*} 
combining the hypothesis and the continuity of $T_a$ in $L^2(\R)$ it follows  that $\theta_2 T_a(\theta_1f)\in L^2(\R))$. Also
\begin{equation}
\label{1}
\begin{aligned}
&\theta_2(x)\,T_a((1-\theta_1)f)(x)\\
&=\int_{-\infty}^{\infty}\theta_2(x)a(x,\xi)\widehat{((1-\theta_1)f)}(\xi)e^{2\pi ix\xi}d\xi\\
&=\int \underset{b(x,\xi)}{ \underbrace{(\int \theta_2(x)a(x,\xi_1+\xi_2)\widehat{(1-\theta_1)}(\xi_1)e^{2\pi i x\xi_1}d\xi_1)}}\widehat{f}(\xi_2)e^{2\pi i x\xi_2}d\xi_2\\
&=T_bf(x)=\int\theta_2(x) k(x,x-z)(1-\theta_1(z))f(z)dz\\
&=\int\theta_2(x) k(x,x-z)(1-\theta_1(z))J^{2m} J^{-2m}f(z)dz
\end{aligned}
\end{equation}
with $-2m<s,\,m\in\Z^+$ and $k(\cdot,\cdot)$ as in \eqref{estimate1}, so integration by parts and Theorem A3 yield the desired result.

\begin{proposition}
\label {P1}
Let $\,f\in L^2(\R)$ and 
\begin{equation*}
J^sf=(1-\partial_x^2)^{s/2}f\in L^2(\{x>a\}) \;\;\;\;\;\;s>0, 
\end{equation*}
then for any $\epsilon>0$ and any $r\in(0,s]$ 
\begin{equation}
\label{inter}
J^{r}f\in L^2(\{x>a+\epsilon\}).
\end{equation}
\end{proposition} 

\noindent{\bf Proof of Proposition \ref{P1}.}
Define
\begin{equation*}
g=J^sf\in L^2(\{x>a\}),
\end{equation*}
thus $\,J^sf\in H^{-s}(\R)$. Let  $\,\theta_j\in C^{\infty}(\R),\,j=1,2,$ with $\,\theta_1(x)=1$ for $x\geq a+\epsilon/4,$ $\,supp\, \theta_1\subset \{x>a\}$, and 
$\,\theta_2(x)=1$ for $x\geq a+\epsilon$ and $\,supp\, \theta_2\subset \{x>a+\epsilon/2\}$,  therefore
$\theta_1 g\in L^2(\R)$. Let $T=J^{i\beta},\,\beta\in\R$. By Lemma \ref{L1} 
\begin{equation*}
\theta_2Tfg=\theta_2J^{s+i \beta}f\in L^2(\R),
\end{equation*}
and since $\,f\in L^2(\R)$ 
\begin{equation*}
\theta_2J^{i\beta}f \in L^2(\R).
\end{equation*}
Hence, by the Three Lines Theorem it follows that 
\begin{equation*}
\theta_2J^{z}f \in L^2(\R),\;\;\;\;\;\;\;z=r+i\beta,\;\;\,r\in[0,s],\;\;\beta\in\R,
\end{equation*}
which completes the proof.

\vskip.1in

\begin{remark} 
\label{R1} In a similar manner one has: for  $\,\epsilon>0$ let $\varphi_{\epsilon}\in C^{\infty}(\R)$ with $\,\varphi_{\epsilon}(x)=1,\;x\geq \epsilon$, $\,supp\,\varphi_{\epsilon}\subset\{x>\epsilon/2\}$ and $\,\varphi_{\epsilon}'(x)\geq 0$. Then
\begin{enumerate}
\item[(I)]  If $\,m\in \Z^+$ and $\varphi_{\epsilon}J^mf\in L^2(\R)$, then  $\,\forall\, \epsilon'>2\epsilon$ 
$$
\varphi_{\epsilon'}\partial_x^jf\in L^2(\R),\;\;\;\;\;\;\;j=0,1,...., m.
$$
\item[(II)]  If $\,m\in \Z^+$ and $\varphi_{\epsilon}\partial_x^jf\in L^2(\R),\;j=0,1,...., m$, then  $\,\forall\,\epsilon'>2\epsilon$ 
$$\varphi_{\epsilon'}J^mf\in L^2(\R).$$

\item[(III)] If $\,s>0,$ and $J^s(\varphi_{\epsilon}\,f)$, $f\in L^2(\R)$, then  $\,\forall\,\epsilon'>2 \epsilon$ 
$$\varphi_{\epsilon'}J^sf\in L^2(\R).$$

\item[(IV)] If $\,s >0,$ and $\varphi_{\epsilon} J^sf$, $f \in L^2(\R)$, then  $\,\forall\,\epsilon'>2\epsilon$  
$$\,J^s(\varphi_{\epsilon'}f)\in L^2(\R).$$

\end{enumerate}
The same results hold with $\,\theta_1,\;\theta_2$, as in \eqref{teta}, instead of  $\chi_{\epsilon},\,\chi_{\epsilon'}$. 
\end{remark}

\vskip.1in
Next, we recall some inequalities obtained in   \cite{KP} :

\begin{TA4}[\cite{KP}]\label{A4}

If $\,s>0$ and $\,p\in(1,\infty)$, then  
\begin{equation}
\label{ine2a}
\|\,J^s(fg)\|_p\leq c(\| f\|_{\infty}\|J^{s} g\|_p + \|J^s f\|_p \|g\|_{\infty}),
\end{equation}
and
\begin{equation}
\label{ine2b}
\begin{aligned}
\|\,[J^s;f]g\|_p&=\| J^s(fg)-fJ^sg\|_p\\
&\leq c(\|\partial f\|_{\infty}\|J^{s-1} g\|_p + \|J^s f\|_p \|g\|_{\infty}).
\end{aligned}
\end{equation}

\end{TA4}
\vskip.1in

Also we shall use the following elementary estimate whose proof is similar to that found in \cite{Fo}, Chapter 6.

\begin{L2}

Let $\,\phi\in C^{\infty}(\R) $ with $\,\phi'\in C_0^{\infty}(\R) $. If $\,s\in \R$, then for any $\,l>|s-1|+1/2$
\begin{equation}
\label{ine3}
\|\,[J^s;\phi]f\|_2 + \|\,[J^{s-1};\phi]\partial_x f\|_2\leq c \,\|J^l\phi'\|_2\,  \|J^{s-1} f\|_2.
\end{equation}
\vskip.2in
\end{L2}

\section{Proof of Theorem \ref{A3}}

 Without loss of generality $\,x_0=0$.   For $\epsilon>0$ and $b\ge 5\epsilon$ define the families of functions
\begin{equation*}
\car,\;\phi_{\epsilon,b},\,\widetilde{\phi_{\epsilon,b}},\,\psi_{\epsilon}\in C^{\infty}(\R),
\end{equation*}
with $\,\car'\ge0$,  $\,\car(x)=0,\;\,x\le \epsilon$, $\,\car(x)=1, \;x\ge b,$
\begin{equation}
\label{000}
\begin{aligned}
\car'(x)&\geq \frac{1}{10(b-\epsilon)} \,1_{[2\epsilon, b-2\epsilon]}(x),\\
\\
& supp(\psi_{\epsilon,b}),\;supp(\widetilde{\psi_{\epsilon,b}}) \subset [\epsilon/4, b],\\
\\
&\phi_{\epsilon,b}(x)=\widetilde{\phi_{\epsilon,b}}(x)=1,\;x\in[\epsilon/2,\epsilon],\\
\\
&\,supp(\psi_{\epsilon})\subset (-\infty,\epsilon/2]\\
\\
&\chi_{\epsilon,b}(x)+ \phi_{\epsilon,b}(x) + \psi_{\epsilon} (x) =1,\;\;\;\;\;\;\;\;\;x\in \R,\\
\\
&\chi^2_{\epsilon,b}(x)+ \widetilde{\phi_{\epsilon,b}}^2(x) + \psi_{\epsilon} (x) =1,\;\;\;\;\;\;\;\;\;x\in \R.
\end{aligned}
\end{equation}
Hence,
$$\;distance (\,supp (\car);\,supp(\psi_{\epsilon}))\geq \epsilon/2.
$$

Formally, we apply the operator $\,J^s\,$  to the equation in \eqref{AA} and multiply by $\,J^s u\chi^2_{\epsilon}(x+vt)$ to obtain 
after integration by parts the identity
\begin{equation}\label{energy}
\begin{split}
&\frac12\frac{d}{dt}\int (J^su)^2(x,t)\chi^2(x+vt)\,dx\\
\\&-\underset{A_1}{ \underbrace{ v\int  (J^s u)^2(x,t)\chi\,\chi'(x+vt)\,dx}}\\
&+\frac32\int  (\partial_x\,J^s u)^2(x,t)\chi \,\chi'(x+vt)\,dx\\
\\ &
-\underset{A_2}{ \underbrace{\frac12\int  (J^s u)^2(x,t)\partial_x^3(\chi^2(x+vt))\,dx}}\\
&+ \underset{A_3}{ \underbrace{\int J^s(u\partial_x u)\,J^su(x,t)\chi^2(x+vt)\,dx}}=0
\end{split}
\end{equation}
where in $\chi$ the index $\epsilon,\,b$ were omitted. We shall do that now on.
\vskip.1in

Case: \underline{$s\in(3/4,1)$.}
\vskip.1in

First, we observe that combining \eqref{notes-2} and the results in section 2 it follows that for any $\,R>0$
\begin{equation}
\label{1a1}
\int^T_0\int_{-R}^R \,|J^ru(x,t)|^2\,dxdt<\infty \;\;\;\;\;\;\;\;\forall \,r\in[0,7/4^+].
\end{equation}
Thus, after integration in  time the terms $A_1$ and $A_2$ in \eqref{energy} are bounded. So it only remains to handle $A_3$.

Thus,
\begin{equation}
\label{z1}
\begin{aligned}
J^s&(u\partial_xu)\chi = J^s(u\partial_xu \chi)-[J^s;\chi](u\partial_xu) \\
&=u\chi J^s\partial_xu +[J^s; u\chi]\partial_xu-[J^s;\chi](u\partial_xu) \\
&=u\chi J^s\partial_xu +[J^s;u\chi]\partial_x(u(\chi+ \phi + \psi))-[J^s;\chi](u\partial_xu)\\
&=B_1+B_2+B_3+B_4+B_5.
\end{aligned}
\end{equation}

Inserting $B_{1} $  in \eqref{energy}  one obtains a term which can be estimated by integration by parts, Gronwall's inequality and \eqref{notes-2}. Using \eqref{ine2b} it follows that
\begin{equation}
\label{z2}
\|B_2\|_2\leq c \|\partial_x(u\chi)\|_{\infty} \| J^s(u\chi)\|_2,
\end{equation}
and
\begin{equation}
\label{z3}
\|B_3\|_2\leq c (\|\partial_x(u\chi)\|_{\infty} \| J^s(u\phi)\|_2+\|\partial_x(u\phi)\|_{\infty} \| J^s(u\chi)\|_2).
\end{equation}
To bound $B_4$ and $B_5$ we apply Corollary \ref{c1} and  \eqref{ine3}, respectively to get
\begin{equation}
\label{z4}
\|B_4\|_2=\| u\chi J^s\partial_x(u\psi)\|_2\leq c \| u\|_{\infty} \|u\|_2
\end{equation}
and
\begin{equation}
\label{z5}
\|B_5\|_2\leq c\|u\|_{\infty} \| u\|_2.
\end{equation}

Collecting the above information \eqref{z1}-\eqref{z5} in \eqref{energy} we obtain \eqref{anotes-4} for any $r\in(3/4,1),\;v>0$ and $\epsilon>0$,
and that for any $v>0,\;\epsilon>0$,
$$
\int_0^T\,\int_{\epsilon-vt}^{R-vt}\,(J^s\partial_xu)^2dxdt<\infty,
$$
from which using the results, Remark \ref{R1}, one obtains \eqref{anotes-5}.

\vskip.1in

Case: \underline{$s\in(m,m+1),\;m\in Z^+$.} 
\vskip.1in

We assume \eqref{anotes-4} and \eqref{anotes-5}  with $s\leq m$. Hence, from the results in section it follows that for any $\,\epsilon>0,\,R>0$ and $\,r\in[0,m]$ 
\begin{equation}
\label {zz1}
\int_0^T\,\int_{\epsilon-vt}^{R-vt}\,(J^r\partial_xu)^2dxdt<\infty.
\end{equation}

 Again the starting point is the energy estimate identity \eqref{energy}. After integrating in time the terms $A_1$ and $A_2$ can be easily bounded using \eqref{zz1}. So it suffices to consider $A_3$. Thus, using the notation introduced in \eqref{000} we have 
\begin{equation}
\label{e1}
\begin{aligned}
\chi & J^s(u\partial_xu)= J^s(u\chi\partial_xu) -\frac{1}{2}[J^s;\chi]\partial_x (u^2)\\
&=u\chi J^s\partial_xu +[J^s;u\chi]\partial_xu -\frac{1}{2}[J^s;\chi]\partial_x (u^2)\\
&=u\chi J^s\partial_xu +[J^s;u\chi]\partial_x(u(\chi+ \phi + \psi))\\
&\;\; -\frac{1}{2}[J^s;\chi]\partial_x ((u^2)(\chi^2+ (\widetilde{\phi})^2 + \psi))\\
&=E_1+E_2+E_3+E_4+E_5+E_6+E_7.
\end{aligned}
\end{equation}

Inserting $E_{1} $ in \eqref{energy} one obtains a term which can be estimated by integration by parts, Gronwall's inequality and \eqref{notes-2}.
From  \eqref{ine2b} we see that
\begin{equation}
\label{e2}
\|E_2\|_2\leq c \|\partial_x(u\chi)\|_{\infty} \| J^s(u\chi)\|_2
\end{equation}
and
\begin{equation}
\label{e3}
\|E_3\|_2\leq c (\|\partial_x(u\chi)\|_{\infty} \| J^s(u\phi)\|_2+\|\partial_x(u\phi)\|_{\infty} \| J^s(u\chi)\|_2).
\end{equation}
For $E_4$  it follows that from Corollary \ref{c1} that
\begin{equation}
\label{e4}
\|E_4\|_2=\| u\chi J^s\partial_x(u\psi)\|_2\leq c \| u\|_{\infty} \|u\|_2.
\end{equation}
  
  For $E_5$ and $E_6$ a combination of \eqref{ine2a} and \eqref{ine3} yields the estimates
\begin{equation}
\label{e5}
\begin{aligned}
\|E_5\|_2&\leq \| [J^s;\chi]\partial_x((u\chi)^2)\|_2\\
&\leq c \|J^s((u\chi)^2)\|_2\leq c\|u\|_{\infty} \| J^s(u\chi)\|_2,
\end{aligned}
\end{equation}
and 
\begin{equation}
\label{e6}
\begin{aligned}
\|E_6\|_2\leq \| [J^s;\chi]\partial_x((u\widetilde{\phi})^2)\|_2=\|J^s((u\widetilde{\phi})^2)\|_2\\
\leq c\|u\|_{\infty} \| J^s(u\widetilde{\phi})\|_2.
\end{aligned}
\end{equation}

Finally, using Corollary \ref{c1} we write
\begin{equation}
\label{e7}
\|E_7\|_2\leq  \| [J^s;\chi]\partial_x((u^2{\psi}))\|_2=\|\chi J^s\partial_x(u^2\psi)\|_2\leq c\|u\|_{\infty}\|u\|_2.
\end{equation}

To complete the estimates in \eqref{e2}, \eqref{e3}, \eqref{e5} and \eqref{e6} we observe that
$$
J^s(u\chi)=J^su \chi +[J^s;\chi](u(\chi+\phi+\psi))=G_1+G_2,
$$
where $G_1$ is the term whose $L^2$-norm we are estimating and $G_2$ is of lower order, (hence bounded by assumption), and
$\,\| J^2(u\phi)\|_2\,$ is bounded by \eqref{anotes-5} (assumption).

Collecting the above information in \eqref{energy} we obtain the desired result.

\vskip.1in


\vspace{5mm}
\noindent\underline{\bf Acknowledgements.} 
\vspace{3mm}

C.E.K.  was supported by the NSF grant DMS-1265249, F.L. was supported by CNPq and FAPERJ/Brazil, and L.V. was supported by an
ERCEA Advanced Grant 2014 669689 - HADE, and by the MINECO
project MTM2014-53850-P.  The authors would like to thank an anonymous referee whose comments helped to improve the presentation of this work.


\end{document}